\documentclass[sigconf, review=false]{acmart}
\renewcommand\footnotetextcopyrightpermission[1]{} 
\renewcommand{\Re}{{\rm I}\!  {\rm R}}
\usepackage{booktabs,algorithm,algcompatible} 
\usepackage{subfigure}
\usepackage[english]{babel}
\usepackage{etoolbox}
\makeatletter
\patchcmd{\maketitle}{\@copyrightspace}{}{}{}
\makeatother

\usepackage[font=small,skip=0pt]{caption}
\setcopyright{none}

\makeatletter
\def\@copyrightspace{\relax}
\makeatother

\begin{document}
\title{A Novel Posture Positioning Method for Multi-Joint Manipulators}


\author{Zhi-Qiang Yao}
\affiliation{%
  \institution{Intelligent Navigation and Remote Sensing Research Centre, Changsha Technology Research Institute of Beidou Industry Safety}
  \streetaddress{Hunan,China}
  \city{Xiangtan,China}
  \postcode{411105}
}
\email{yaozhiqiang@xtu.edu.cn}
\renewcommand{\shortauthors}{Karimpanal}

\author{Yi-Jue Dai}
\affiliation{%
  \institution{The Chinese University of Hong Kong}
  \streetaddress{Shenzhen, China}
  \city{Shenzhen,China}
  \postcode{518127}
}
\email{218019041@link.cuhk.edu.cn}
\renewcommand{\shortauthors}{Karimpanal}

\author{Qing-Na Li}
\affiliation{%
  \institution{School of Mathematics and Statistics/Beijing Key Laboratory on MCSSCI}
  \streetaddress{Beijing,China}
  \city{Beijing,China}
  \postcode{100081}
}
\email{qnl@bit.edu.cn}
\renewcommand{\shortauthors}{Karimpanal}

\author{Dang Xie}
\affiliation{%
  \institution{Intelligent Navigation and Remote Sensing Research Centre}
  \streetaddress{Hunan,China}
  \city{Xiangtan,China}
  \postcode{411105}
}
\email{249689394@qq.com}
\renewcommand{\shortauthors}{Karimpanal}

\author{Zehui Liu}
\affiliation{%
  \institution{Intelligent Navigation and Remote Sensing Research Centre}
  \streetaddress{Hunan,China}
  \city{Xiangtan,China}
  \postcode{411105}
}
\email{740662158@qq.com}
\renewcommand{\shortauthors}{Karimpanal}


\begin{abstract}
Safety and automatic control are extremely important when operating manipulators. For large engineering manipulators, the main challenge is to accurately recognize the posture of all arm segments. In classical sensing methods, the accuracy of an inclinometer is easily affected by the elastic deformation   in the manipulator's arms. This results in  big error accumulations  when sensing the angle of joints between arms one by one. In addition, the sensing method based on machine vision is not suitable for such kind of outdoor working situation yet. In this paper, we propose a novel posture positioning method for multi-joint manipulators based on wireless sensor network localization.  The posture sensing problem is formulated as a Nearest-Euclidean-Distance-Matrix (NEDM) model. The resulting approach is referred to as EDM-based posture positioning approach (EPP) and it  satisfies the following guiding principles: (i) The posture of each arm segment on a multi-joint manipulator must be estimated as accurately as possible; (ii) The approach must be computationally fast; (iii) The designed approach should not be susceptible to obstructions.
To further improve accuracy, we explore the inherent structure of manipulators, i.e., fixed-arm length. This is naturally presented as linear constraints in the NEDM model. For concrete pumps, a typical multi-joint manipulator, the mechanical property that all arm segments always lie in a 2D plane is used for dimension-reduction operation. Simulation and experimental results show that the proposed method provides   efficient solutions for posture sensing problem and can obtain preferable localization performance with faster speed than applying the existing localization methods.
\end{abstract}
%
%

\keywords{Multi-joint Manipulators, Concrete Pumps, Sensor Network Localization, Euclidean Distance Matrices}

\copyrightyear{2018}
\acmYear{2018}
\setcopyright{none}
\acmConference[GECCO '18 Companion]{Genetic and Evolutionary
Computation Conference Companion}{July 15--19, 2018}{Kyoto, Japan}
\acmBooktitle{GECCO '18 Companion: Genetic and Evolutionary
Computation Conference Companion, July 15--19, 2018, Kyoto,
Japan}\acmDOI{10.1145/3205651.3208762}
\acmISBN{978-1-4503-5764-7/18/07}

\maketitle

\section{Introduction}
\label{sec:introduction}
Multi-joint manipulator is a large-scale arm-like mechanism that consists of several arm segments.
Safety and automatic control are extremely important when operating such manipulators. For large engineering manipulators, the main challenge is to accurately recognize the posture of all arm segments without collisions.
Take the concrete pump as an example. Due to large errors in operating systems, people need to drag the hose (on the end of the arm) to the desired position. This is usually unstable and may lead to serious accidents \cite{2007Zhou-928-932}.  Collisions often occur   between arms of two different manipulators, or between an arm of a manipulator and other objects (such as people, trees and buildings) \cite{2017WYang}. Therefore, the posture of manipulators needs to be accurately recognized. This is referred as posture sensing problem.

Existing methods for the posture sensing problem can be summarized into two types. The first type is based on angles of rotation. The angles are usually measured by inclination sensors   on every arm joint  \cite{20101, 20102, 2008Wang-494-497}.  However, the accuracy of an inclinometer is easily affected by the elastic deformation   in the manipulator's arms. It results in  big error accumulations when sensing the angle of joints between arms one by one. 
The other type relies on machine vision, which is the technology  to provide imaging-based automatic inspection. For example, Li et al. \cite{2010G} scanned QR codes on joints through cameras to recognize postures. Mila \cite{2010P} and Vorobieva et al. \cite{2010H} extracted contour of objects through an expensive binocular vision system. Such type of methods is not suitable for outdoor working situations, since there may be obstacles in the operating environment \cite{2016sjia}.


In this paper, we propose a novel posture positioning method for multi-joint manipulators based on wireless sensor network localization (SNL). The posture sensing problem is formulated as a Nearest-Euclidean-Distance-Matrix (NEDM) model  and the resulting approach is referred to as EDM-based posture positioning approach (EPP). The proposed approach  satisfies the following guiding principles: (i) The posture of each arm segment on a multi-joint manipulator must be estimated as accurately as possible; (ii) The approach must be computationally fast; (iii) The   approach should not be susceptible to obstructions. To further improve accuracy, we explore the inherent structure of manipulators, i.e., fixed-arm length. This is naturally presented as linear constraints in the NEDM model. For concrete pumps, a typical multi-joint manipulator, the mechanical property that all arm segments always lie in a 2D plane is used for dimension-reduction operation. Simulation and experimental results show that the proposed method provides  efficient solutions for posture sensing problem and can obtain preferable localization performance with less cputime  than applying the existing localization methods.


The rest of the paper is organized as follows. In Section \uppercase\expandafter{\romannumeral2}, the mathematical statement for the posture sensing problem is introduced, followed by the NEDM model. 
 In Section \uppercase\expandafter{\romannumeral3}, we study  the concrete pump as an example. A coordinate transformation technique is introduced to tackle the special feature of concrete pumps. We conduct extensive simulations in Section \uppercase\expandafter{\romannumeral4}, and draw some conclusions in Section \uppercase\expandafter{\romannumeral5}.

\section{Mathematical Model for Multi-Joint Localizations}
\label{methods}
In this section, we will give the mathematical statement  and NEDM model for the posture sensing problem.

 Let $\mathcal S^n$ denote the space of $n$ by $n$ symmetric matrices. A positive semidefinite matrix $B\in\mathcal S^n$ is denoted as $B\succeq 0$. We use $\|\cdot\|$ to denote the $l_2$ norm for vectors, and the Frobenius norm for matrices. The small (capital) variables  are vectors (matrices). Let $\hbox{Diag}(x)$ denote the diagonal matrix formed by vector $x\in \Re^n$, and  $\hbox{diag}(X)$ denote the vector formed by the diagonal elements in matrix $X\in \mathcal S^n$.  Let $\textbf{e}\in\Re^n$ be the vector with all elements one.

{\bf Problem Statement.} The posture sensing problem is to recognize the posture of each arm segment on manipulators. An efficient way is to find the position of each joint. From this point of view, it can be reformulated as a wireless sensor network localization problem. Assume that there are $p$ arm segments on a manipulator, denoted as the $1$-st segment, $2$-nd segment, $\dots$, $p$-th segment from the end to the final base. The end of the first segment is denoted as point $1$, and the joint between segment $1$ and segment $2$ is denoted as point $2$, etc. The point where the $p$-th arm segment joined with the turntable is denoted as point $p+1$. The aim of posture sensing is to locate the position of point $1, 2, \dots, p$.

The mathematical statement of posture sensing problem is   described as follows. Given the locations ${\textbf{w}}_{p+1}, \dots, {\textbf{w}}_n\in\Re^3$, and estimated locations $\tilde {\textbf{w}}_1, \dots, \tilde {\textbf{w}}_p$ at time $t-1$, try to find the locations ${\textbf{w}}_1, \dots, {\textbf{w}}_p\in\Re^3$ of $1, \dots, p$ at the current time $t$.  
The available information also includes some noised distances $\delta_{ij}$, $i\in\{1, \dots, p\}$, $j\in\{p+1, \dots, n\}$. See Fig.\ref{fig:2} for the demonstration of the problem (take $p=5,n=9$ as example).
Moreover,  the posture sensing problem also enjoys the following inherent features: the length of each arm segment is fixed and  is available. In other words,
\begin{equation}\label{feature-1}
\hbox{\textbf{Feature I: }}\ \|{\textbf{w}}_{j+1}-{\textbf{w}}_j\|=L_j, \ j = 1, \dots, p,
\end{equation}
where $L_j$ is the length of the $j$-th arm segment, $j=1, \dots, p$.
 \begin{figure}[htpb]
 	\begin{center}		
    \includegraphics[width=3.5in]{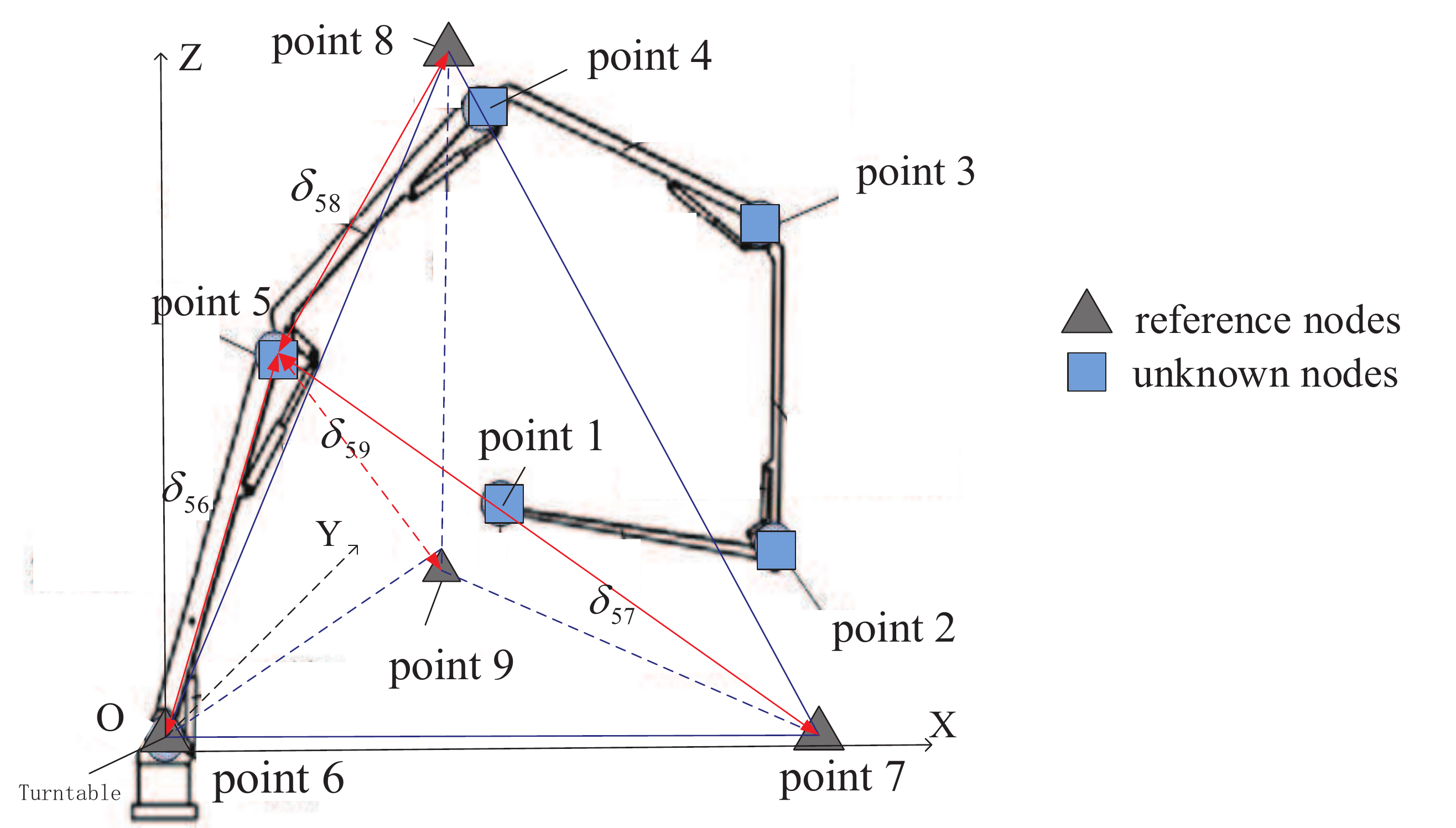}
	\caption{The demonstration of the posture sensing problem of multi-joint manipulators.\label{fig:2}}
	\label{fig:2}
	\end{center}		
 \end{figure}


{\bf An NEDM Model.} To propose our NEDM model, assume that
 the posture changes between any  successive test time are small.
 A Euclidean Distance Matrix (EDM) $ {D}\in\mathcal S^n$ for a set of points $\{{\textbf{w}}_1, \dots, {\textbf{w}}_n\}\subseteq \Re^m$ is defined as the squared distance matrix, which is $D_{ij}=\|{\textbf{w}}_i-{\textbf{w}}_j\|^2$, $i,j = 1, \dots, n$. For  posture sensing problem, we are looking for an EDM  $D$, which is generated by $\textbf{w}_1,\dots, \textbf{w}_n\in\Re^3$. 
 Putting it in another way, suppose we have a noised EDM matrix $ {G}\in\mathcal S^n$,  we are looking for an EDM with corresponding points in  $\Re^3$, which is nearest to $\ {G}$. This can be formulated as a least-squares problem
 \begin{equation}
\begin{aligned}
\min\limits_{\mathrm{ {D}}\in \mathcal{S}^{n}}
& \frac{1}{2}\|\mathrm{ {H}}\circ (\mathrm{ {D}}-\mathrm{ {G}})\|_F^2:=f(D)\\
\hbox{s.t.}~&\hbox{diag}( {D})= {0}, \ ~ {D}\in K^n_-,\\
&D_{jj+1}=L_j^2,~j=1,\dots, p,\\
&rank(JDJ)\leq 3.
\end{aligned}
\label{equ:mod}
\end{equation}
where "$\circ$" denotes the Hadamard product, and $ {H}\in\mathcal S^n$ is a prescribed weight matrix. 
 The   constraints in the first line describe that $D$ is an EDM, where  $K^n_{-}$ is defined by $K^n_-=\{ {X}\in \mathcal S^n|~ \textbf{v}^TX\textbf{ v}\geq 0, ~\forall\ \textbf{v}^T \textbf{e}=0\}$ (See  \cite{Schoenberg35, Young1938, Qi} for more details). 
 The second line in the constraints tackles the fix-length arm feature. The third line of the constraints guarantees the points that generate the EDM $D$ are in $\Re^3$. 

After obtaining the matrix $ {D}$, the embedding points $\textbf{v}_1, \dots, \textbf{v}_n$ can be given by classical multidimensional scaling (cMDS), which is given below.

First, a spectral decomposition is conducted as follows
\begin{equation}\label{cMDS1}
-\frac12 {JDJ} =  {P}\hbox{Diag}(\lambda_1,\dots, \lambda_n) {P}^T,  \  {J} =  {I}-\frac1n \textbf{ee}^T,
\end{equation}
then
\begin{equation}\label{cMDS2}
 [\textbf{v}_1 \ \textbf{v}_2\ \dots \ \textbf{v}_n] ^T =  {P}_1\hbox{Diag}(\lambda_1^\frac12,\dots,\lambda_r^\frac12) ~\in\Re^{n\times r},
\end{equation}
where $r=3$ is the required dimension, 
 $ {I}$ is the identity matrix of $n$ by $n$, $\lambda_1\ge\lambda_2\ge\dots\ge\lambda_n\ge0$ are eigenvalues and $ {P}_1\in\Re^{n\times r}$ consists of  the corresponding eigenvectors of $\lambda_1, \ldots, \lambda_r$.

The post process is then the same as that for SNL. The final estimations of $1,\dots, p$ are denoted as $\hat{ \textbf{w}}_1,\dots, \hat {\textbf{w}}_p$. See \cite{Toh} for more details of the post process.

 {\bf Remark 1.} Here we would like to emphasize that the advantage of the NEDM model is that it can efficiently deal with linear constraints in low computational cost. Problem (\ref{equ:mod}) is a nonconvex optimization problem due to the rank constraint. Solutions to (\ref{equ:mod}) are recently studied by Qi, et al in \cite{2014Qi-351-351} and \cite{Qi}, where an efficient majorized penalty method has been proposed. See \cite{2014Qi-351-351, Qi} for more details.

{\bf Remark 2.}  Compared with other existing SNL solvers, such as the SDP approach in \cite{2006Biswas-360-371} as well as SR-LS in \cite{2008Beck-1770-1778}, our NEDM model fully make use of the fixed-arm length feature, which will further improve the accuracy of localization. This will be verified later in numerical part. 

{{\bf Properties of NEDM Model.}
The majorized penalty method  for NEDM model (\ref{equ:mod}) is to solve the following    penalty problem
 \begin{equation}
\begin{aligned}
\min\limits_{\mathrm{ {D}}\in \mathcal{S}^{n}}
&  f(D)+ cq(D)\\
\hbox{s.t.}~&\hbox{diag}( {D})= {0}, \ ~ {D}\in K^n_-,\\
&D_{jj+1}=L_j^2,~j=1,\dots, p.\\
\end{aligned}
\label{equ:mod-major}
\end{equation}
where $c>0$ is the penalty parameter, $q(D)$ is the penalty function for nonconvex rank constraint $rank { {J}} { {D}} { {J}})\leq 3$. We refer to \cite{2014Qi-351-351} for more discussions on the choices of $q(D)$. We have the following result addressing the EDM solutions of (\ref{equ:mod}) and (\ref{equ:mod-major}), which comes from \cite[Prop. 5.1]{2014Qi-351-351} and originally summarized from Prop. 3.1 and Prop. 3.2 in \cite{Gao-Sun}.

{\bf Proposition 1.}
Let $D^*_c$ denote a global optimal EDM solution of (\ref{equ:mod-major}). Let $D_r$ be a feasible EDM solution of (\ref{equ:mod}) and $D^*$ be an optimal EDM solution of the following convex problem
 \begin{equation}
\begin{aligned}
\min\limits_{ { {D}}\in \mathcal{S}^{n}}
&  f(D) \\
\hbox{s.t.}~&\hbox{diag}( {D})= {0}, \ ~ {D}\in K^n_-\\
&D_{jj+1}=L_j^2,~j=1,\dots, p.\\
\end{aligned}
\label{equ:mod-convex}
\end{equation}
\begin{itemize}
\item [(i)] If $\hbox{rank}(JDJ)\le 3$, 
then EDM $D^*_c$ already solves (\ref{equ:mod}).  That is, it is the optimal EDM solution of (\ref{equ:mod}).
\item [(ii)] If $c$ is chosen to satisfy $c\ge (f(D_r)-f(D^*))/\epsilon$ for some given $\epsilon>0$, then we have
\[
p(D^*_c)\le \epsilon, \hbox{ and } f(D^*_c)\le \nu^*-cp(D^*_c),
\]
where $\nu^*$ is the optimal function value of (\ref{equ:mod}).

\end{itemize}


The result in (ii) means that when the rank error measured by $q(\cdot)$ at $D^*_c$ is less than $\epsilon$, 
EDM $D^*_c$ is an $\epsilon$-optimal EDM solution \cite{Gao-Sun} of (\ref{equ:mod}).

Define another weight matrix $W:=\hbox{Diag(u)}$, where $u\in\Re^{n}$ is defined by
\[
u_i=\hbox{max}\{\tau, \max\{H_{ij}, \ j = 1,\ldots, n\}\}, \ i = 1, \ldots, n.
\]
for some $\tau>0$.

The following result address the convergence property of the majorized penalty method, which comes from Proposition 5.3 in \cite{2014Qi-351-351}.

{\bf Proposition 2.}
Let $\{D^k\}$ be the sequence generated by the majorized penalty method. Then $\{f(D^k)\}$ is a monotonically decreasing sequence. If $D^{k+1}=D^k$ for some $D^k$, then $D^k$ is an optimal EDM solution of (\ref{equ:mod}). Otherwise, the infinite sequence $\{D^k\}$ satisfies
\[
\frac12\|W^{\frac12}(D^{k+1}-D^k)W^{\frac12}\|^2\le f(D^k)-f(D^{k+1}), \ k =0, 1, \ldots
\]
Moreover, the sequence $\{D^k\}$ is bounded and any accumulation point is a $B$-stationary point of (\ref{equ:mod-major}).

 Proposition 2 basically says that the proposed method either terminates at the optimal EDM solution of (\ref{equ:mod}) or  any of the accumulation points is a B-stationary point of the penalty problem (\ref{equ:mod-major}). Note that a B-stationary point is usually the best point that a numerical method can find for (\ref{equ:mod-major}) as it is nonconvex.
 }
%

\section{An Example of Concrete Pumps}
A concrete pump, one of the largest engineering manipulators, is a device to transfer liquid concrete through pumping. It is widely used in construction projects. A common structure of the concrete pump is that it has only one turntable, meaning that the horizontal angles of rotation of all arm segments are the same, and the major motions are stretching and folding. In other words, besides feature  \eqref{feature-1}, concrete pumps also enjoy the following feature
\[
\hbox{\textbf{Feature II:} All arm segments lie in the same vertical plane.}
\]

Mathematically speaking,  points $1,\dots, p,\ p+1 $ should lie in a 2D plane. However, the estimated points   returned by the NEDM model  (\ref{equ:mod})  may not guarantee to be in a 2D plane. In the following, we propose a new approach based on the NEDM model. The idea is as follows. 
To make use of feature II, we first map all the data onto the vertical plane via coordinate transformation, then apply the NEDM model (\ref{equ:mod}) in 2D space. The resulting method is referred to as Coordinate-EDM-based Posture Positioning approach (CEPP). 
We give the details below.

\subsection{Step 1: Coordinate Transformation}
\begin{figure}[htpb]
    \begin{center}
    \includegraphics[width=2.8in]{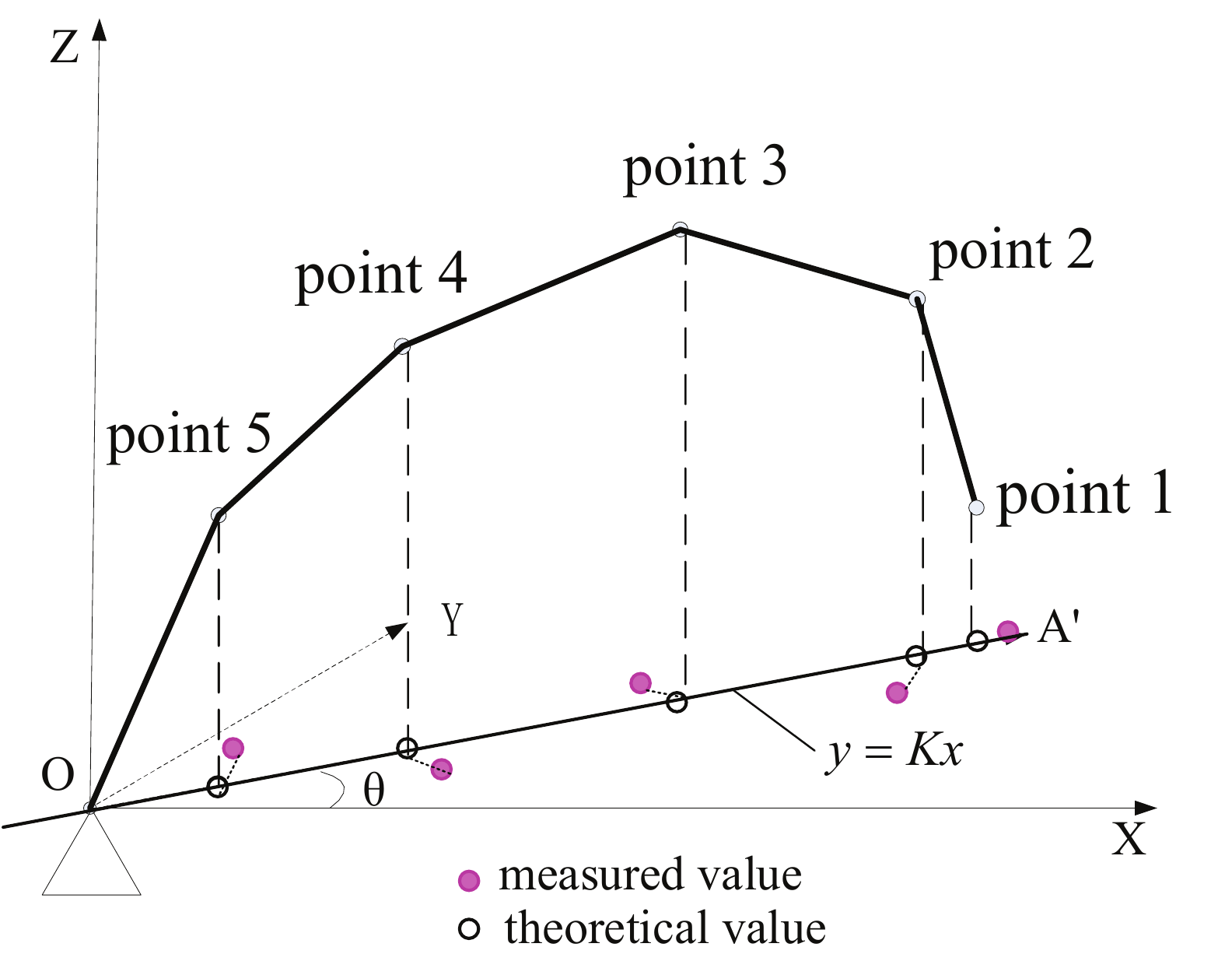}
    \end{center}
    \caption{Determine the vertical plane of the concrete pump.}
    \label{fig:3}
\end{figure}

 We use a coordinate transformation technique to map the points onto the vertical plane where the $p$ arm segments lie in. To facilitate our statement, we put the origin at point $p+1$, and denote it as $O$. Recall that all points lie in $\Re^3$, and $\tilde {\textbf{w}}_1,\dots, \tilde {\textbf{w}}_p$ are the estimations of $1,\dots, p$ at time $t-1$. The way to determine the vertical plane is shown in Fig. \ref{fig:3}. We first map $\tilde {\textbf{w}}_1, \dots, \tilde {\textbf{w}}_p$ onto the horizontal plane $XOY$. Ideally, the projections will lie in a line. However, due to the estimated errors, this may fail. We use the least squares fitting to find a line $OA'$  whose slope is denoted by  $  K $.  Then the angle of rotation between the vertical plane $A'OZ$ and $XOZ$ is $\theta = \hbox{arc}\hbox{tan} K$. The vertical plane $A'OZ$ is the 2D plane that we are looking for.

\begin{figure}[htpb]
	\begin{center}
		\includegraphics[width=3in]{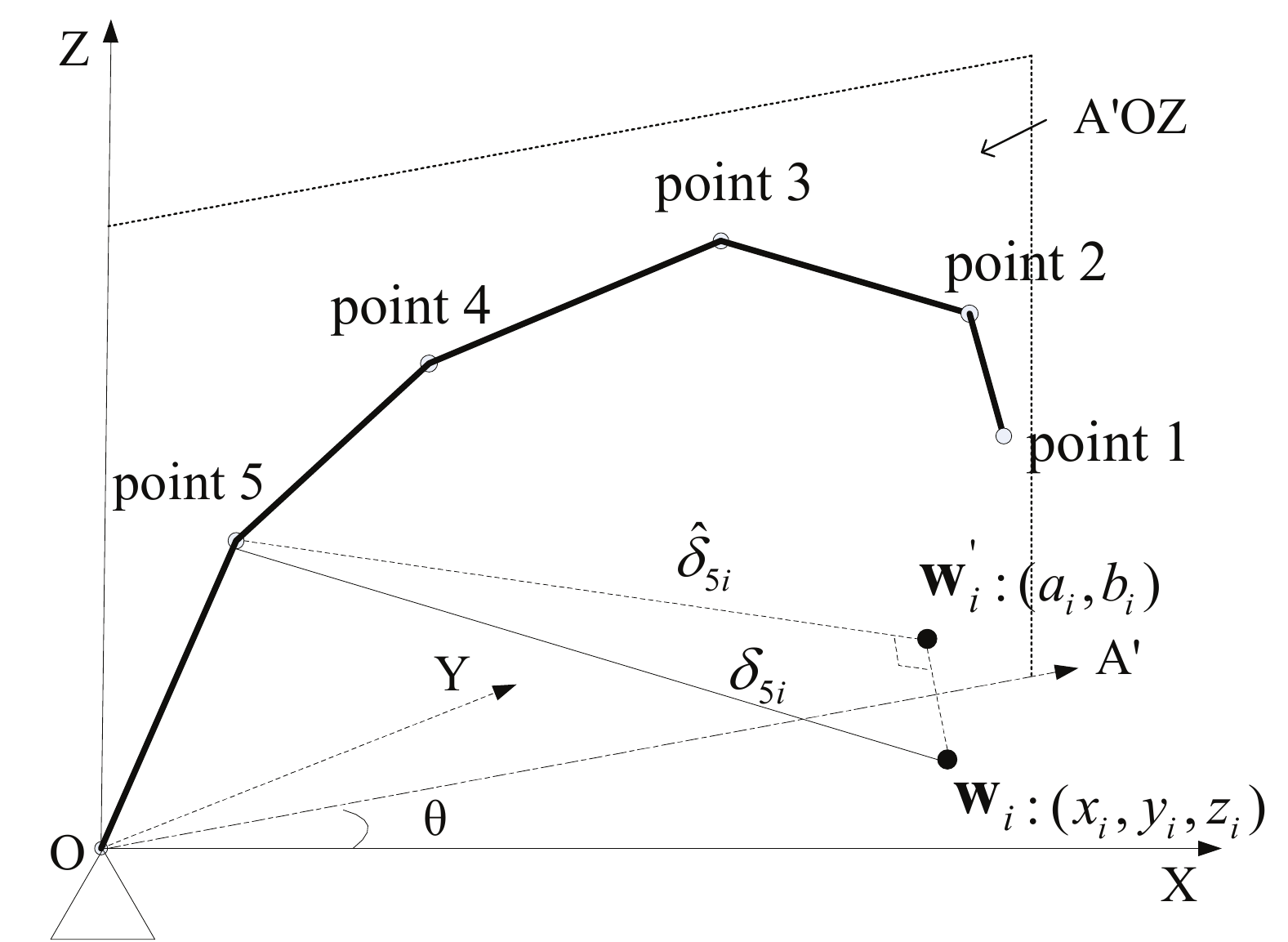}
	\end{center}
	\caption{The proposed coordinate transformation of point $i$ and the noised distance $\delta_{ij}$.}
	\label{fig:4}
\end{figure}

Next, we project points $\tilde {\textbf{w}}_1,\dots, \tilde {\textbf{w}}_p$ and $  {\textbf{w}}_{p+1},\dots,   {\textbf{w}}_{n}$ onto the vertical plane $A'OZ$ to get the projection points in 2D plane, denoted as  $\tilde {\textbf{w}}_1',\dots, \tilde {\textbf{w}}_p'\in\Re^2$ and $  {\textbf{w}}_{p+1}',\dots,   {\textbf{w}}_{n}'\in \Re^2$.  For any  point  $\textbf{w}=(x ,y , z )\in \Re^3$, it can be mapped onto $A'OZ$ by the following coordinate transformation formulae to get new points $\textbf w'=(a ,b )\in\Re^2$
 \begin{equation}
\left\{
\begin{aligned}
&{{a }} = \sqrt {x ^2 + ({x {\tan \theta}})^2} + \sin \theta ({y } - {x }\tan \theta )  \\
&{{b } = {z }},
\end{aligned}
\right.
\label{equ:trans1}
\end{equation}
The observed noisy distances $\delta_{ij}$ can also be transferred to the  corresponding distances $\hat \delta_{ij}$ in $A'OZ$ by $$\hat\delta_{ij} =\sqrt {\delta_{ij}^2 - x_i^2 - y_i^2 + a_i^2},\ i = 1,\dots, p, \ j = p+1, \dots,n.$$
The coordinate transformation and the distance projection are illustrated in Fig.\ref{fig:4}.

\subsection{Step 2: Dimension Reduction via NEDM Model}
  Having obtained the points and distances in the vertical plane $A'OZ$, we can apply the NEDM model in 2D space by solving problem (\ref{equ:mod}) with
 the rank constraint  replaced by $\hbox{rank}(\textbf{JDJ})\le2$. After the post process, we get estimates $\hat {\textbf{w}}_1', \dots, \hat {\textbf{w}}_p'\in\Re^2$ and can be easily transformed back to 3D positions, denoted as $\hat {\textbf{w}}_1, \dots, \hat {\textbf{w}}_p\in\Re^3$.

 {\bf Remark 3.} Here, we would like to highlight that our   approach belongs to SNL approaches and therefore is independent of the inclination angle and machine vision. As a result, the proposed approach can get rid of the error accumulation caused by inclination angles. Moreover, since we only use the location information, there is no computer vision involved. It is particularly suitable for large manipulators such as concrete pumps.

 \section{Simulation Results}
In this section, we test our approach and compare with other methods via simulations. This section is divided into three parts. In the first part, we do simulations on a large manipulator. In the second part, we test the case of concrete pumps. In the last part, we provide semi-physical test to further verify the efficiency of the proposed method.

\subsection{Large Manipulator Posture Recognition}
In this part, we test our approach on a large manipulator with five arm segments $(p=5)$, the lengths of arm segments from turntable are set to 9m, 7m, 7m, 9m, 9m, respectively. {We would like to emphasize that due to the large errors for angle-rotation based approaches, as well as the inapplicability of computer vision, we only compare different solvers of SNL approaches.}
They are detailed as follows. To see the role of the equality constraints, we test  NEDM model (\ref{equ:mod}) without the fixed-arm constraints. The resulting approach is denoted as EPP1. We use EPP2 to denotes the NEDM model (\ref{equ:mod}), which takes account of the equality constraints\footnote{ {Available from  http://www.personal.soton.ac.uk/hdqi/,  last access date: } {   July 1, 2020.}}. We also compare with some classical solvers for SNL, including the SDP approach\footnote{Available from http://web.stanford.edu/~yyye/Col.html,  last access date:   July 1, 2020.}  \cite{2006Biswas-360-371} and SR-LS approach\footnote{ Available from https://github.com/daiyijue-XTU/EMBED-SRLS, last access date:   July 1, 2020.} \cite{2008Beck-1770-1778}. Note that both the SDP approach and SR-LS approach do not consider the fixed-arm length constraints. For EPP1 and EPP2, one particular choice of $ {G}$ is as follows.
\begin{equation}
\mathrm{ {G}}=	\left\{
\begin{aligned}
&\|\mathrm{\textbf{w}}_i-\mathrm{\textbf{w}}_j\|^2;\quad i,j=p+1,\dots, n,\\
&\delta_{ij}^2;\quad i=1,\dots, p, j=p+1, \dots, n,\\
&\|\tilde{\mathrm{ \textbf{w}}}_i-\tilde{\mathrm{\textbf{w}}}_j\|^2;\quad i,j=1,\dots, p.
\end{aligned}
\right.
\label{equ:EDM}
\end{equation}

\begin{figure}[t!]
 	\begin{center}			
 		\includegraphics[height=2.5in,width=3.5in]{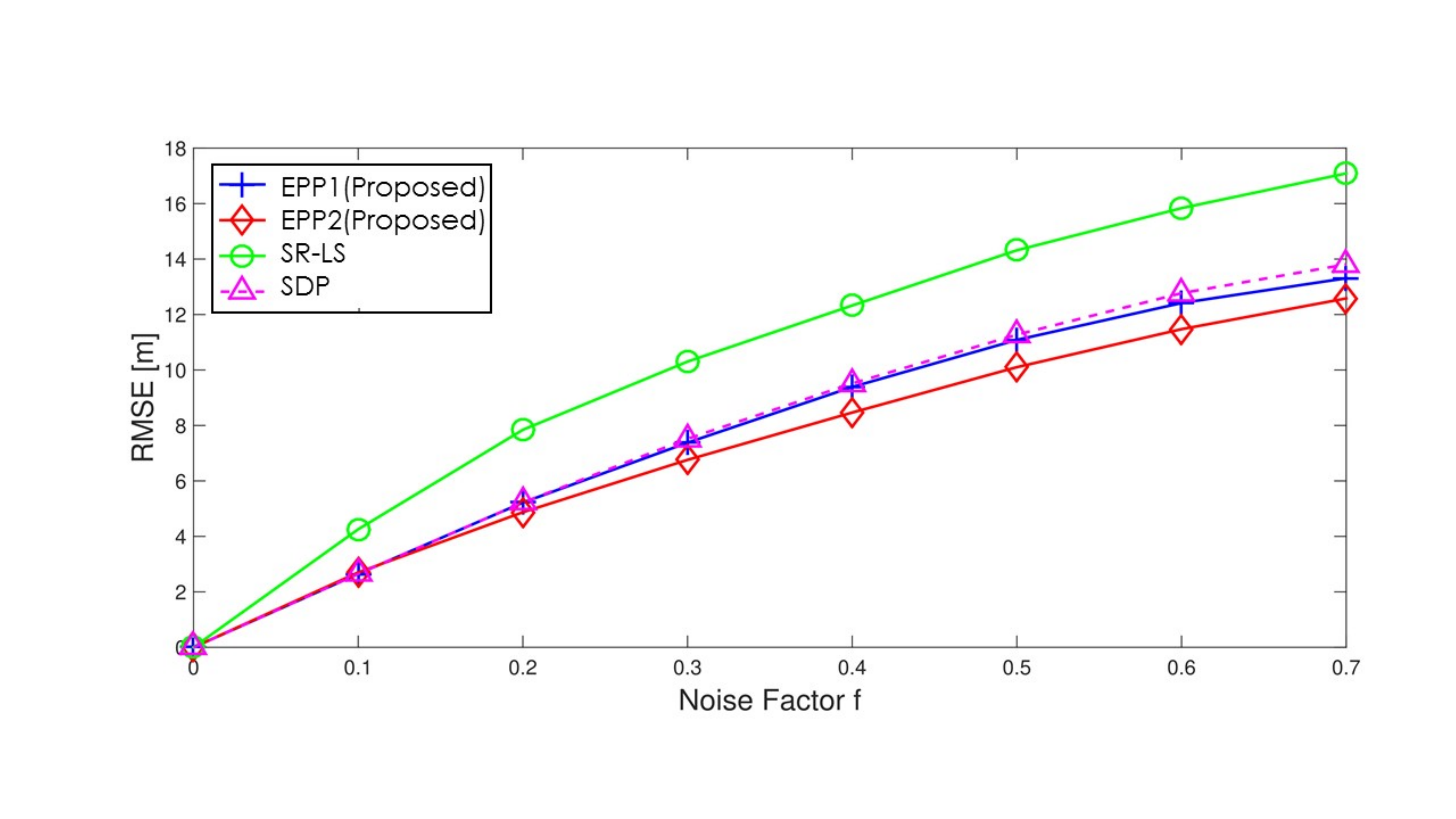}
 	\end{center}		
 	\caption{The comparison of approaches on RMSE with different noise factor.}
 	\label{fig:dd}
 \end{figure}

 The matrix $ {H}$ is taken to be the matrix with all entries equal to one. Random error is added to an observed distance in the following manner:  $$\delta_{ij}=\|\textbf{w}_i-\textbf{w}_j\|\cdot(1+N(0,1)*\eta),$$ where noise factor $\eta$ is a given number between [0,1] and $N(0,1)$ denotes a standard normal random variable.
   The RMSE values over 3000 runs with different $\eta$ are shown in Fig.\ref{fig:dd}, where  $$RMSE := 1/\sqrt{p}(\sum_{i=1}^{p}\|\hat{\textbf{w}}_i-{\textbf{w}}_i\|^2)^{1/2.}.$$ One can see that EPP1 gives quite similar RMSE as SDP does, and both are better than those given by SR-LS. Comparing EPP2 with EPP1, the performance has been improved after adding the equality constraints. The average CPU time and absolute error of the length of every arm segment are presented in Table \ref{tab:tab1}. Obviously, in a very low time consumption, EPP2 works well and outperforms other methods, both in RMSE and absolute error.

  \begin{table}[t!]
 	\centering
 	\caption{Comparison of different approaches on the absolute error of arm length  [\upshape{m}]}\label{tab:tab1}	
 	\begin{tabular}{ccccccc}		
 		\toprule
 		SEGMENT &1 & 2 & 3 & 4 & 5 & CPU time [s]\\
 		\midrule
 		SDP&1.639&1.623&1.554&1.540&1.681&0.286\\
 		\midrule
 		SR-LS&2.215&4.964&5.780&4.743&5.024&0.519\\
 		\midrule
 		EPP1&1.630&1.629&1.554&1.548&1.692&0.041\\
 		\midrule
 		\textbf{EPP2}&\textbf{0.057}&\textbf{0.014}&\textbf{0.015}&\textbf{0.016}&\textbf{0.017}&\textbf{0.047}\\
 		\bottomrule	
 	\end{tabular}			
 \end{table}

\subsection{Concrete Pumps Posture Positioning}
In this part, we compare the proposed method to the traditional posture sensing method that are based on angles of rotation. In \cite{20101,20102,2008Wang-494-497}, the posture is recognized through trigonometric operations based on the rotational angles between two adjacent arms. Here we use TPSM to represent this traditional posture sensing method, and CEPP2 to present the method that we apply the coordinate transformation before EPP2.

\begin{figure}[t!]
	\begin{center}			
		\includegraphics[height=2.5in,width=3.5in]{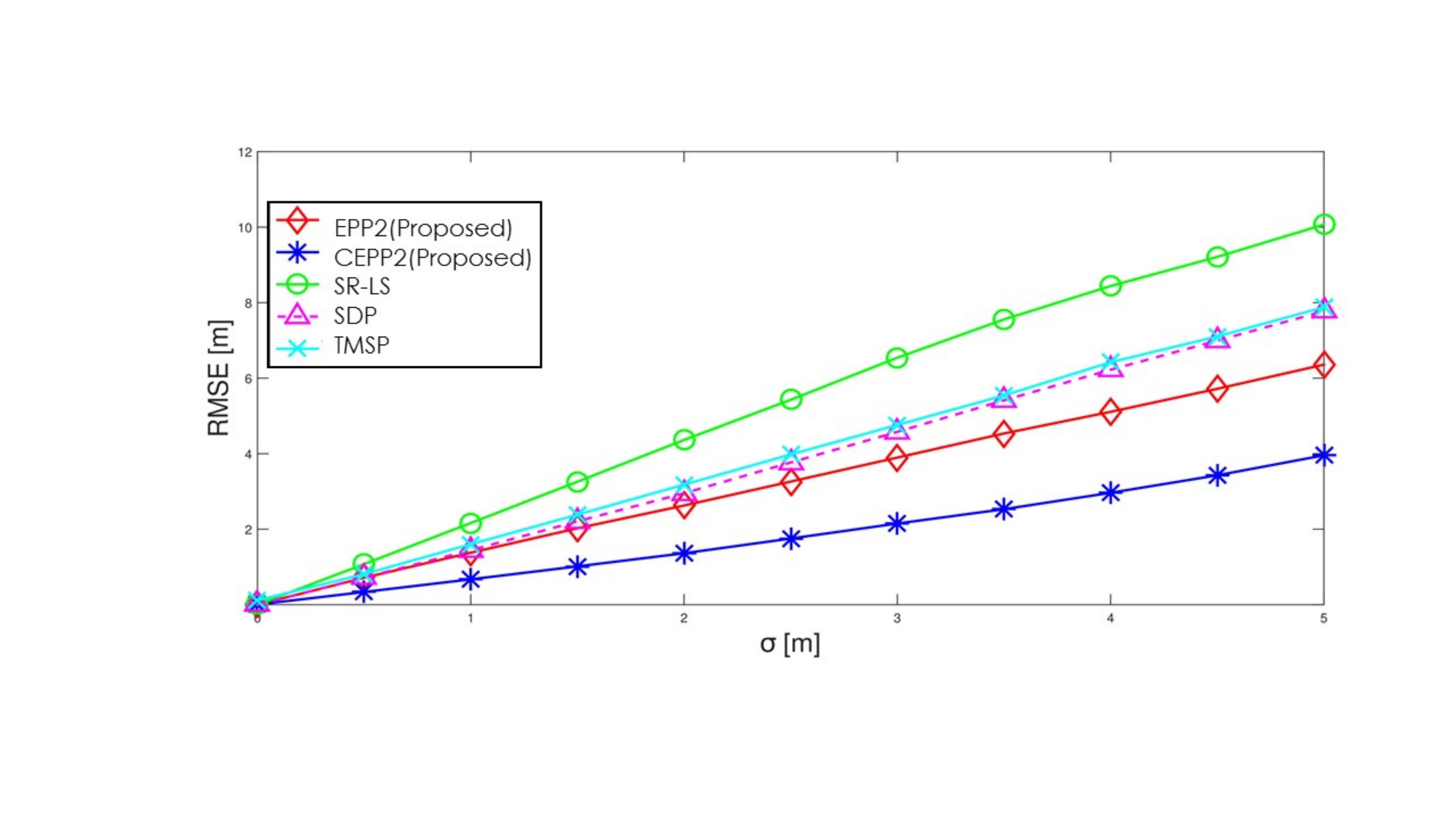}
	\end{center}		
	\caption{The comparison on RMSE with different standard deviation value.}
	\label{fig:6}
\end{figure}

\begin{figure}[]
	\begin{center}			
		\includegraphics[height=2.5in,width=3.5in]{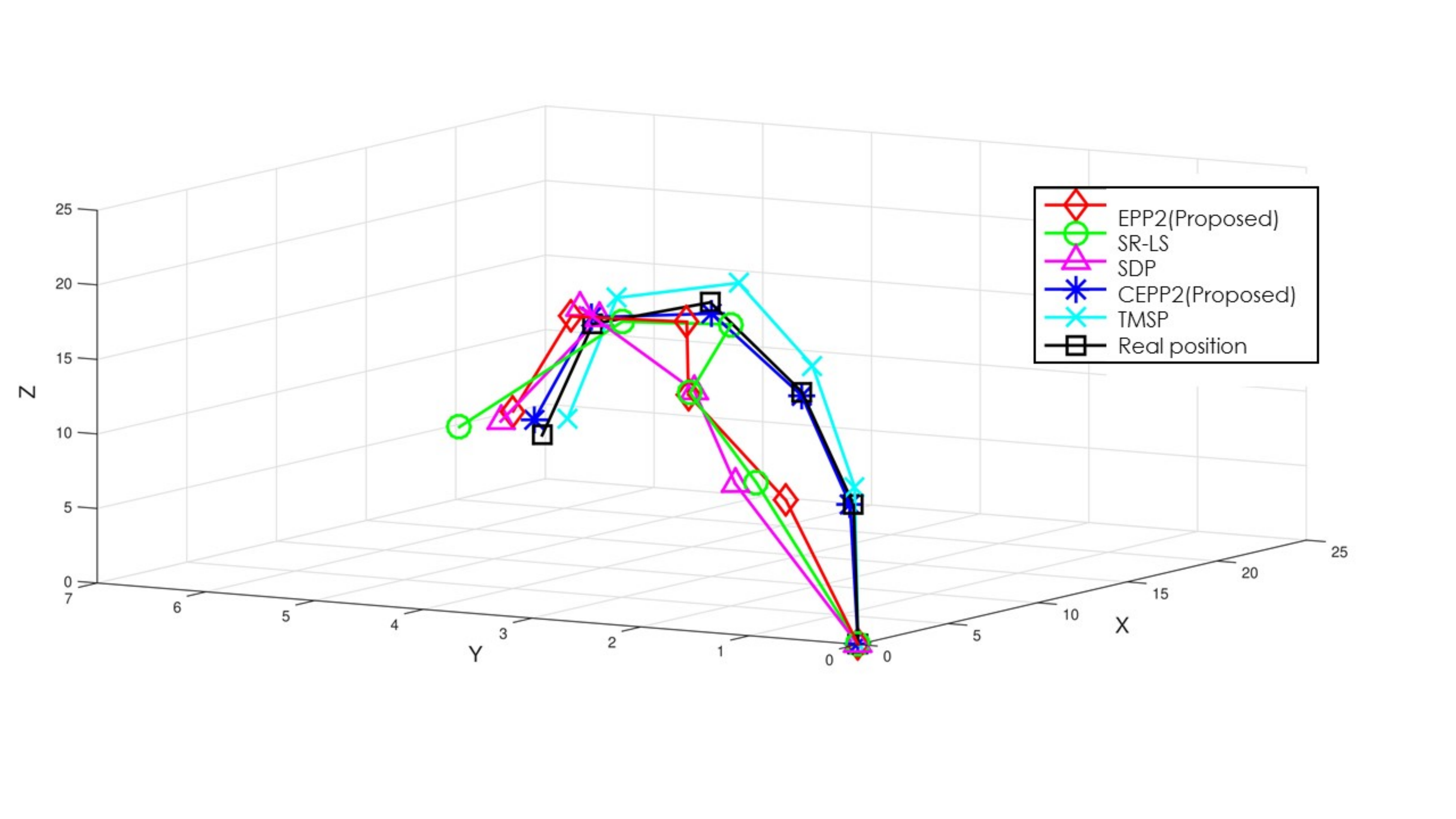}
	\end{center}		
	\caption{Real posture vs. estimated postures from different methods $(\sigma=2)$.}
	\label{fig:7}
\end{figure}

{\bf TEST One}. In this test, the additive noise term in the measurements follow i.i.d $N(0,\sigma^2)$. A comparison on RMSE with different standard deviation $\sigma$ over 3000 runs is shown in Fig.\ref{fig:6}, where the angle error of the turntable is set to be $\pm 0.5^\circ$ in TPSM, and other rotation errors are transferred to the corresponding distance error components. It shows that the coordinate transformation  is helpful to the process, and the performance of CEPP2 is superior to EPP2 and other estimates for all values of $\sigma$. The postures of the concrete pump estimated from different approaches are shown in Fig.\ref{fig:7}. Apparently, the posture of CEPP2 obviously outperforms other methods (see the line marked with rhombus).

{\bf TEST Two}. Assume there exist some NLOS (Non-Line Of Sight) paths between points $p+1,\dots,n$ and points $1,\dots,p$. We add some NLOS error terms to the measured data generated in the first test. It is modeled as an exponential random variable with parameter $\gamma$, as in the examples of  \cite{2017Qi-187-187}. Let $\sigma=0.2m$ and $\gamma=2m$. The comparison on RMSE and average CPU time are shown in Table \ref{tab:tab3}. The results show that the CEPP2   achieved   significantly better performance than other methods.

 \begin{table}[h]
  \centering
		\caption{Comparison of different methods in NLOS scenario.}	\label{tab:tab3}
		\begin{tabular}{cccccc}		
			\toprule
			& SDP & SR-LS & TPSM & EPP2 & \textbf{CEPP2}\\
			\midrule
			RMSE [m]&4.503&5.464&4.765&4.249&\textbf{2.533}\\
			\midrule
			CPU time [s]& 0.312&0.539&0.019&0.048&\textbf{0.063}\\		
			\bottomrule	
		\end{tabular}			
\end{table}

\subsection{Semi-Physical Experiment for Concrete Pumps}

\begin{figure}[h]
	\centering
	\subfigure[]
	{ 	
		\includegraphics[height=4cm,width=4cm]{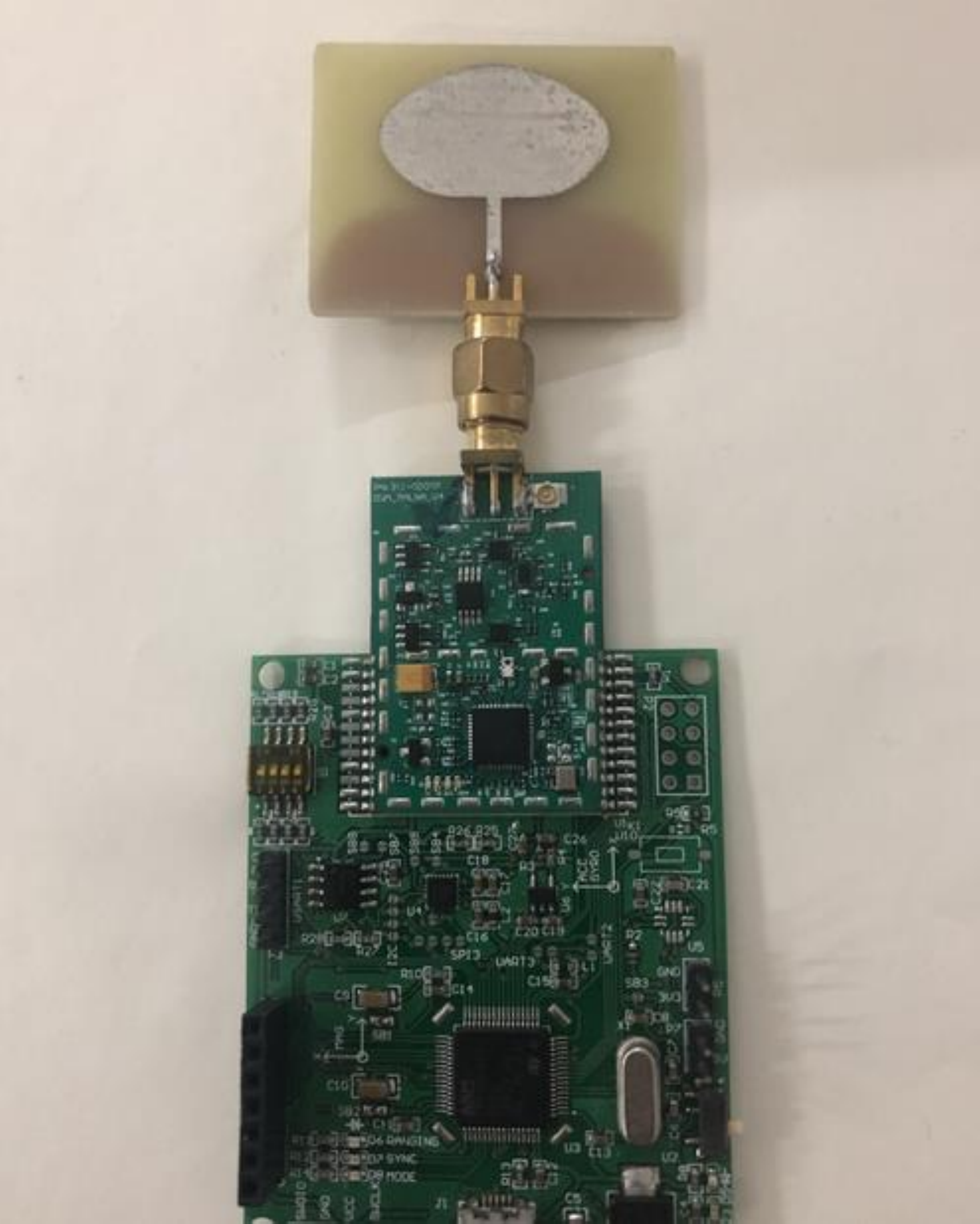}
		
	}
	\subfigure[]
	{ 	
		\includegraphics[height=4cm, width=4cm]{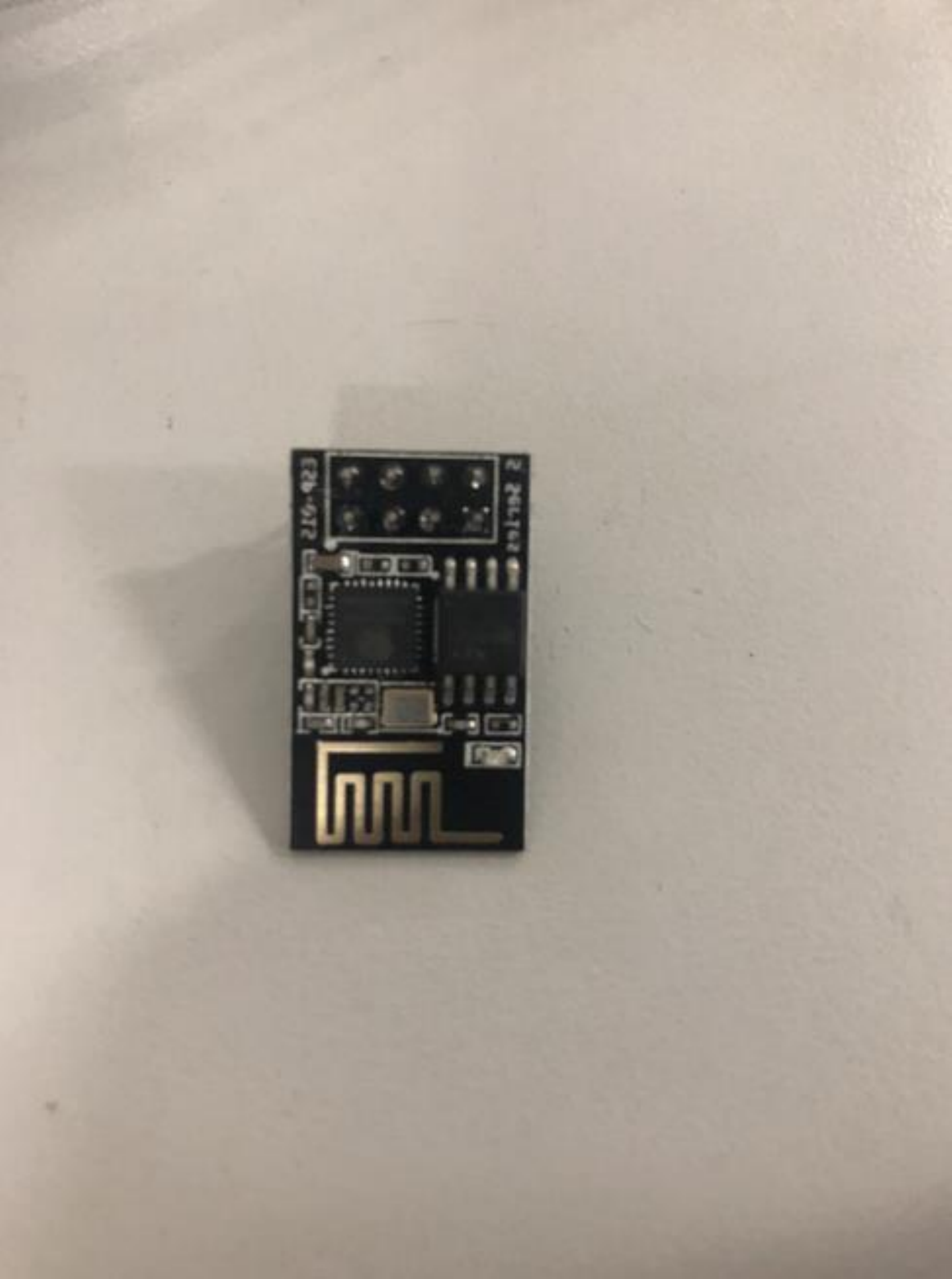}

	}
	\caption{Physical drawing of UWB: dwm1000 and WiFi: ESP8226}
	\label{fig:12}
\end{figure}

\begin{figure}[h]
	\centering
	{ 	
		\includegraphics[width=7cm]{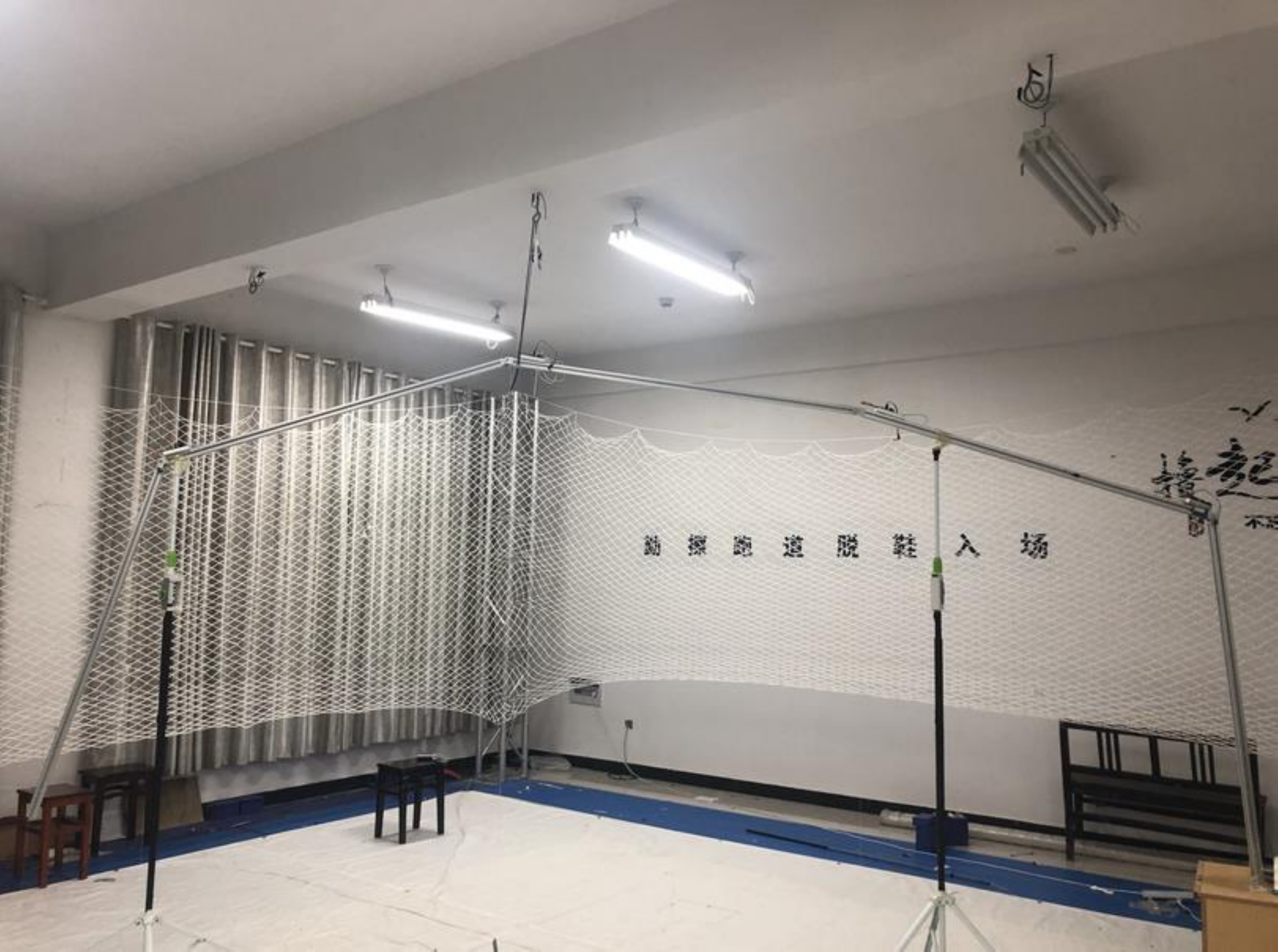}
		
	}
	{ 	
		\includegraphics[width=7cm]{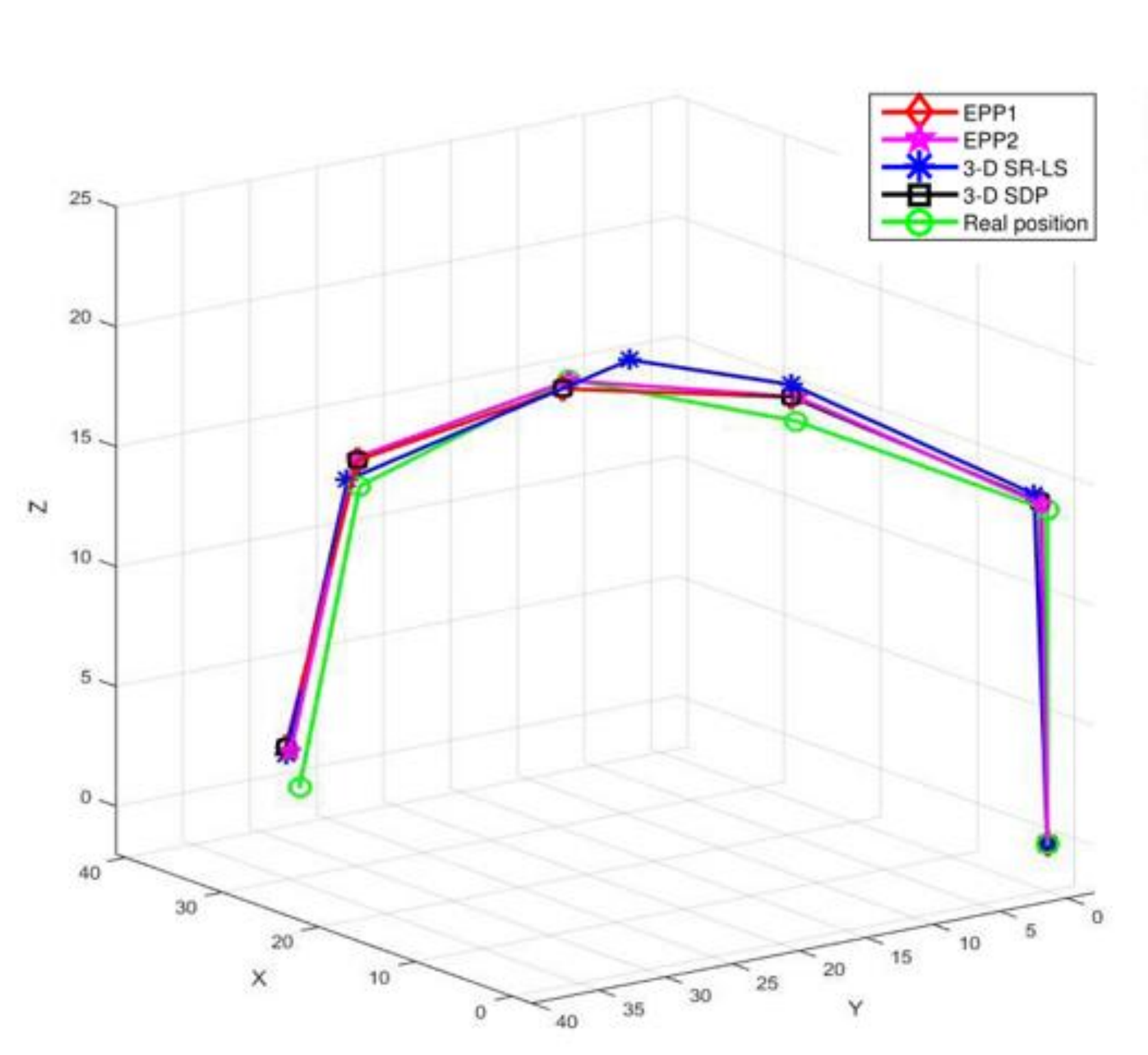}
		\label{fig:11}
	}
	\caption{Experimental scenarios of the test on concrete pumps and the results}
	\label{fig:10}
\end{figure}

In order to test our proposed methods in real data, we build a semi-physical concrete pump by scaling down at the same proportion of the real industrial concrete pump with  nodes as shown in the top figure of Fig.\ref{fig:10}.  To measure the distances, we use Ultra-Wide Bandwidth (UWB) modules as shown in Fig. \ref{fig:12} (a).  We use the WIFI as shown in Fig. \ref{fig:12} (b)  as  the wireless communication equipment.  The test is conducted as follows.  As shown in  Fig. \ref{fig:2}, we used nine UWB modules, each of which was placed  on the point 1 to 9. We  used  ten WIFI modules  for wireless  communication.  The UWB units placed on points 1 to point 5 (5 targets) represents the locations which  need to be localized, and they are in the same vertical plane. The rest   4  units represents points 6 to point 9 with known positions  (4 anchors). Then we use the range method to measure the distances among each pair of  points, and transmit the location information to   computer through wifi module shown in Fig.\ref{fig:12} (b). In our test, we set the coordinate for the 4 anchors from point 6 to 9 are [0, 0, 0], [0, 420, 0], [400, 420, 0], [221, 167, 247] (cm), respectively. The measured distance matrix is as follows (unit: cm)

\begin{equation*}
G = \begin{bmatrix}
0      &       \\
387 & 0&  \\
377  & 498& 0\\
381 &  609&  347&0\\
244 &  539&  375&172&0\\
100 &  440& 365 &284&141&0\\
49 &  329& 381&396&257&131&0\\
169 & 207 & 351 &433&352&243&137&0\\
305 & 109 &  390&521&457&367&246&141&0\\
\end{bmatrix}
\end{equation*}
where the matrix G is symmetric, and $G_{ij}$ denotes the squared distance between the point $i$ and point $j$. The arm lengths of our experiment concrete pump are all 150cm. The comparison of different SNL approaches on the absolute error of arm lengths are shown in Table \ref{tab:tab2}. The average RMSE of various algorithms to the positioning results are listed in Table \ref{tab:tab4}. Sicne CEPP2 applies coordinate transformation based on EPP2. Therefore CEPP2 is more accurate than EPP2. In CEPP2, all joints are transformed to the same plane for simplify the calculation.  Therefore, the cputime for  CEPP2  is the smallest. The test results  shows that our proposed NEDM model and coordinate transformation technique   worked quite well in tested environment, compared with other SNL solvers.

  \begin{table}[t!]
	\centering
	\caption{Comparison of different appraoches on the absolute error of arm length [\upshape{cm}]}\label{tab:tab2}	
	\begin{tabular}{ccccccc}		
		\toprule
		SEGMENT &1 & 2 & 3 & 4 & 5 & CPU time [s]\\
		\midrule
		SDP&0.672&0.488&0.603&0.684&0.492&0.296\\
		\midrule
		SR-LS&1.358&1.348&1.191&1.241&1.777&0.583\\
		\midrule
		EPP1&0.673&0.487&0.602&0.684&0.493&0.044\\
		\midrule
		\textbf{EPP2}&\textbf{0.002}&\textbf{0.011}&\textbf{0.008}&\textbf{0.001}&\textbf{0.001}&\textbf{0.047}\\
		\bottomrule	
	\end{tabular}			
\end{table}

  \begin{table}[t!]
	\centering
	\caption{Comparison of different approaches on the average error of joint position  [\upshape{cm}]}	\label{tab:tab4}
	\begin{tabular}{ccccccc}		
		\toprule
		SEGMENT &1 & 2 & 3 & 4 & 5 & average RMSE\\
		\midrule
		SDP&1.11&1.58&1.21&0.96&1.65&1.30\\
		\midrule
		SR-LS&1.50&1.31&1.47&0.68&2.83&1.56\\
		\midrule
		EPP1&2.89&1.10&1.40&1.11&2.19&1.74\\
		\midrule
		\textbf{EPP2}&1.05&1.34&1.21&1.13&1.02&\textbf{1.16}\\
		\bottomrule	
	\end{tabular}			
\end{table}

\subsection{Further Discussions}
{Compared with the indoor environment as tested in Section C, one may wonder whether the performance of the proposed approach can be maintained when applied outdoors.  We address this issue from the following three aspects. Firstly, for normal outdoor environment,   the proposed approach is more favorable in such situations. The reason is that there is serious mulit-path effect of signals in indoor environment, making the noise more complicated than the outdoor situation. The performance of the new approach will be better than that presented in Section C. In other words, more reliable and  accurate postures will be provided via the new approach when applied outdoors. Secondly, even when there are obstacles among the sensors in practice, the approach still works. Because the new approach is designed based on wireless communication, rather than machine vision.  Finally, in some very extreme cases, for example, there is electric welding which is very close to the manipulator,  it will generate severe noises to the wireless signals.  Such unusual noises will bring down the performance of both our approach and the machine vision based approach. 
}

\section{Conclusions}
 In this paper, we solved the posture sensing  problem through the wireless network localization  approach and we introduced the NEDM model   for large engineering manipulators. To further improve the performance, some linear constraints are added in the NEDM model to tackle the inherent feature of fixed arm length. In the case study of the concrete pump, we proposed a   coordinate-EDM  posture positioning approach to further tackle the feature that all joints of arms lies in a 2D space. 
  Simulation and experimental results verified the accuracy and efficiency of the proposed methods. Our work brought an innovative way to solve the posture sensing problem for manipulators and will make it increasingly possible to realize the automatic pouring and safe operation, especially for   concrete pumps.

\end{document}